\documentclass[12pt]{amsart}
\numberwithin{equation}{section}
\newtheorem{theorem}{Theorem}

\newtheorem{proposition}{Proposition}
\numberwithin{theorem}{section}
\numberwithin{lemma}{section}
\numberwithin{proposition}{section}
\newtheorem{corollary}{Corollary}
\newcommand{\ignore}[1]{}
\newcommand{\w}{\omega}
\newcommand{\R}{{\bf R}}
\newcommand{\beq}{\begin{equation}}
\newcommand{\enq}{\end{equation}}

\begin{document}
\tracingpages 1
\title[local regularity]{\bf
When does a Schr\"odinger heat equation permit positive solutions }
\author{Qi S. Zhang}
\address{Department of Mathematics, , University of
California, Riverside, CA 92521, USA }
\date{January 2006}

\begin{abstract}
We introduce some new  classes of time dependent functions whose
defining properties take into account of oscillations around
singularities. We study properties of solutions to the heat equation
with coefficients in these classes which are much more singular than
those allowed under the current theory. In the case of $L^2$
potentials and $L^2$ solutions, we give a characterization of
potentials which allow the Schr\"odinger heat equation to have a
positive solution. This provides a new result on the long running
problem of identifying potentials permitting a positive solution to
the Schr\"odinger equation.

We also establish a nearly necessary and sufficient condition on
certain sign changing potentials such that the corresponding heat
kernel has Gaussian upper and lower bound.

Some applications to the Navier-Stokes equations are given. In
particular, we derive a new type of a priori estimate for solutions
of Navier-Stokes equations. The point is that the gap between this
estimate and a sufficient condition for all time smoothness of the
solution is logarithmic .

\end{abstract}
\maketitle
\tableofcontents
\section{Introduction}

In the first part of the paper we would like to study the heat
equation with a singular $L^2$ potential $V = V(x, t)$ , i.e.
\begin{equation}
\label{1.1}
\begin{cases}
 \Delta u + V u - u_t=0, \quad \text{in} \quad {\bf R}^n \times
(0, \infty), \quad V \in L^2({\bf R}^n \times
(0, \infty)), \ n \ge 3, \\
u(x, 0) = u_0(x), \quad x \in {\bf R}^n, \quad u_0 \in L^2(\R^n).
\end{cases}
\end{equation} Since we are only concerned with local
regularity issue in this paper, we will always assume that $V$ is
zero outside of a  cylinder in space time: $B(0, R_0) \times [0,
T_0]$, unless stated otherwise. Here $R_0$ and $T_0$ are  fixed
positive number. The $L^2$ condition on the potential $V$ is modeled
after the three dimensional vorticity equation derived from the
Navier-Stokes equation. There the potential $V$ is in fact the
gradient of the velocity which is known to be in $L^2$. The unknown
function $u$ in (\ref{1.1}) corresponds to the vorticity which is
also known to be a $L^2$ function.

We will use the following definition of weak solutions.

\medskip

\noindent {\bf Definition 1.1}. Let $T>0$. We say that $u \in
L^1_{loc}(\R^n \times (0,, T))$ is a solution to (\ref{1.1}) if
$V(.) u(., t) \in L^1_{loc}(\R^n \times (0, T))$ and
\[
\int_{\R^n}  u_0(x)  \phi(x)  dx + \int^{T}_0 \int_{\R^n} u \phi_t
dxdt + \int^{T}_0 \int_D  u  \Delta \phi dxdt + \int^T_0 \int_{\R^n}
V u \phi dxdt =0
\]for all smooth, compactly supported $\phi$ vanishing on $\R^n
\times \{ T \}.$

It is well known that $L^2$ potentials in general are too singular
to allow weak solutions of (\ref{1.1}) to be bounded or unique.
Therefor further assumptions must be imposed in order to establish a
regularity theory. The classical condition on the potential $V$ for
H\"older continuity and uniqueness of weak solutions is that $V \in
L^{p, q}_{loc}$ with $\frac{n}{p} + \frac{2}{q} < 2$. This condition
is sharp in general since one can easily construct a counter
example. For instance for $V = a/|x|^2$ with $a>0$, then there is no
bounded positive solution to (\ref{1.1}) (c.f. [BG]). In fact in
that paper, it is shown that if $a$ is sufficiently large, then even
weak positive solutions can not exist. There is a long history of
finding larger class of potentials such that some regularity of the
weak solutions is possible. Among them is the Kato class, time
independent or otherwise. Roughly speaking a function is in a Kato
type class if the convolution of the absolute value of the function
and the fundamental solution of Laplace or the heat equation is
bounded. This class of functions are moderately more general than
the standard $L^{p, q}$ class. However, it is still far from enough
for applications in such places as the vorticity equation mentioned
above. We refer the reader to the papers [AS], [Si], [Z], [LS] and
reference therein for results in this direction. The main results
there is the continuity of weak solutions with potentials in the
Kato class. In addition, equation (\ref{1.1}) with $V$ in Morrey or
Besov classes are also studied. However, these classes are
essentially logarithmic improvements over the standard $L^{p, q}$
class. In the paper [St], K. Sturm proved Gaussian upper and lower
bound for the fundamental solution when the potential belongs to a
class of time independent, singular oscillating functions. His
condition is on the $L^1$ bound of the fundamental solution of a
slightly "larger" potential.

In this paper we introduce a new class of {\it time dependent}
potentials which can be written as a nonlinear combination of
derivatives of a function. The general idea of studying elliptic and
parabolic equations with potentials as the spatial derivative of
some functions is not new. This has been used in the classical books
[LSU], [GT] and [Lieb]. Here we also allow the appearance of time
derivative which can not be dominated by the Laplace operator.
Another innovation  is the use of a suitable combination of
derivatives. The class we are going to define in section 2
essentially characterize all $L^2$ potentials which allow
(\ref{1.1}) to have positive $L^2$ solutions.

The question of whether the Laplace or the heat equation with a
potential possesses a positive solution has been a long standing
one. For the Laplace equation, when the potential has only mild
singularity, i.e. in the Kato class, a satisfactory answer can be
found in the Allegretto-Piepenbrink theory. See Theorem C.8.1 in the
survey paper [Si].
 Brezis and J. L. Lions (see [BG] p122) asked when (\ref{1.1}) with
more singular potential has a positive solutions. This problem was
solved in [BG] when $V = a/|x|^2$ with $a>0$. In the case of general
time independent potentials $V \ge 0$, it was solved in [CM] and
later generalized in [GZ]. However the case of time dependent or
sign changing potentials is completely open. One of the main results
of the paper (Theorem 2.1) gives a solution of the problem with
$L^2$ potentials.  The main advantage of the new class
 of potentials is that it correctly captures the cancelation
effect of sign changing functions. Moreover, we show in Theorem
2.2-3 below that, if we just narrow the class a little, then the
weak fundamental solutions not only exist but also have Gaussian
upper bound. A Gaussian lower bound is also established under
further but necessary restrictions.

Some of the results of the paper can be generalized beyond $L^2$
potentials. However we will not seek full generalization  this time.

Before proceeding further  let us fix some notations and symbols, to
which will refer the reader going over the rest of the paper.

\medskip
\noindent {\bf Notations}. We will use $\R^+$ to denote $(0,
\infty)$. The letter $C$, $c$ with or without index will denote
generic positive constants whose value may change from line to line,
unless specified otherwise.  When we say  a time dependent function
is in $L^2$ we mean its square is integrable in $\R^n \times \R^+$.
We use $G_V$ to denote the fundamental solution of (\ref{1.1}) if it
exists. Please see the next section for its existence and
uniqueness. The symbol $G_0$ will denote the fundamental solution of
the heat equation free of potentials. Give $b>0$, we will use $g_b$
to denote a Gaussian with $b$ as the exponential parameter, i.e.
\[
g_b = g_b(x, t; y, s) = \frac{1}{(t-s)^{n/2}} e^{-b |x-y|^2/(t-s)}.
\]Given a $L^1_{loc}$ function $f$ in space time, we will use $g_b
\star f(x, t)$ to denote
\[
\int^t_0\int_{\R^n}g_b(x, t; y, s) f(y, s) dyds.
\]When we say that $G_V$ has Gaussian upper bound, we mean that
exists $b>0$ and $c>0$ such that $G_V(x, t; y, s) \le c g_b(x, t; y,
s)$. The same goes for the Gaussian lower bound.

When we say a function is a positive solution to (\ref{1.1}) we mean
it is a nonnegative weak solution which is not identically zero.

\medskip

Here is the plan of the paper. In the next section we provide the
definitions, statements and proofs of the main results. In section
3, we define another class of singular potentials, called heat
bounded class. Some applications to the Navier-Stokes equation is
given in Section 4.

\section{singular potentials as combinations of derivatives}

\subsection{Definitions, Statements of Theorems}

$$
$$

\medskip

\noindent {\it {\bf Definition 2.1.} Given two  functions $V \in
L^2(\R^n \times \R^+)$, $f \in L^1_{loc}(\R^n \times \R^+)$ and
$\alpha>0$, we say that
\[
V = \Delta f - \alpha |\nabla f|^2 - f_t,
\]if there exists sequences of functions $\{ V_i \}$ and
$\{ f_i \}$  such that the following conditions hold for all $i=1,
2...$:

\noindent (i).   $V_i \in L^2(\R^n \times \R^+)$, \  $\Delta f_i \in
L^2(\R^n \times \R^+)$, \ $\partial_t f_i \in L^2(\R^n \times
\R^+)$, $f_i \in C(\R^n \times \R^+)$.

\noindent (ii).  $V_i \to V$ strongly in $L^2(\R^n \times \R^+)$,
$|V_i| \le |V_{i+1}|$, $f_i \to f$ a.e. and $f_i(x, 0) = f_{i+1}(x,
0)$.

\noindent (iii).
\[ V_i= \Delta f_i - \alpha |\nabla f_i|^2 -
\partial_t f_i.
\]}
\medskip

Here we remark that we do not assume $\Delta f$, $|\nabla f|^2$ or
$\partial_t f$ are in $L^2$ individually. This explains the lengthy
appearance of the definition.

\medskip
The main results of Section 2 are the next three theorems. The first
one states a necessary and sufficient condition such that
(\ref{1.1}) possesses a positive solution. The second theorem
establishes Gaussian upper and lower bound for the fundamental
solutions of (\ref{1.1}).  The third theorem is an application of
the second one in the more traditional setting of $L^{p, q}$
conditions on the potential. It will show that our conditions are
genuinely much broader than the traditional ones.

It should be made clear that there is no claim on uniqueness in any
of the theorems. In the absence of uniqueness how does one define
the fundamental solution? This is possible due to the uniqueness of
problem (\ref{1.1}) when the potential $V$ is truncated from above.
This fact is proved in Proposition 2.1 below. Consequently we can
state
\medskip

\noindent {\it {\bf Definition 2.2.} The fundamental solution $G_V$
is defined as the pointwise limit of the (increasing) sequence of
the fundamental solution $G_{V_i}$ where $V_i = \min \{ V, \ i \}$
with $i=1, 2, ...$.}
\medskip

We remark that $G_V$ thus defined may be infinity somewhere or
everywhere. However we will show that they have better behavior or
even Gaussian bounds under further conditions.

\begin{theorem}
(i). Suppose (\ref{1.1}) with some $u_0 \ge 0$ has a positive
solution. Then
\[
V = \Delta f - |\nabla f|^2 - \partial_t f,
\]with
 $e^{-f} \in L^2(\R^n \times \R^+)$. Moreover $f \in L^1_{loc}(\R^n \times \R^+)$
 if $\ln u_0 \in L^1_{loc}(\R^n)$.

(ii). Suppose
\[
V = \Delta f - |\nabla f|^2 - \partial_t f
\]for some $f$ such that $e^{-f} \in L^2(\R^n \times \R^+)$. Then
the equation in (\ref{1.1}) has a positive $L^2$ solution for some
$u_0 \in L^2(\R^n)$.
\end{theorem}
\medskip

\medskip
\begin{theorem}
(i). Suppose
\[
V = \Delta f - \alpha |\nabla f|^2 - \partial_t f
\]for one given $\alpha>1$ and $f \in L^\infty(\R^n \times \R^+)$.
Then $G_V$ has Gaussian upper bound in all space time.

(ii). Under the same assumption as in (i),  if $g_{1/4} \star
|\nabla f|^2 \in L^\infty(\R^n \times \R^+)$, then $G_V$ has
Gaussian lower bound in all space time.

(iii). Under the same assumption as in (i), suppose $G_V$ has
Gaussian lower bound in all space time. Then there exists  $b>0$
such that
\[
g_b \star |\nabla f|^2 \in L^\infty(\R^n \times \R^+).
\]
\end{theorem}
\medskip

{\it Remark 2.1.} At the first glance, Theorem 2.1 may seem like a
restatement of existence of positive solutions without much work.
However Theorem 2.2 shows that if one just puts a little more
restriction on the potential $V$, then the fundamental solution
actually has a Gaussian upper bound. Under an additional but
necessary assumption, a Gaussian lower bound also holds. Even  the
widely studied potential  $a/|x|^2$  in $\R^n$ can be recast in the
form of Theorem 1.1, as indicated in the following

{\it Example 2.1.} For a real number $b$, we write $f = b \ln r$
with $r = |x|$. Then direct calculation shows, for $r \neq 0$,
\[
\frac{b(n-2-b)}{r^2} = \Delta f - |\nabla f|^2.
\]Let $a =b(n-2-b)$. Then it is clear that the range of $a$ is $(-\infty,
(n-2)^2/4]$. In this interval (1.1) with $V= a/|x|^2$ permits
positive solutions. This recovers the existence part in the
classical result [BG]. Highly singular, time dependent examples can
be constructed by taking $f=\sin (\frac{1}{|x|-\sqrt{t}})$ e.g.

Moreover the corollaries below relate our class of potentials with
the traditional "form bounded" or domination class (\ref{domi}),
(see also \cite{Si:1} below). In the difficult time dependent case,
Corollary 1 shows that  potentials permitting positive solutions,
can be written as the sum of one form bounded potentials and the
time derivative of a function almost bounded from above by a
constant.
\medskip

{\it Remark 2.2.} From the proof, it will be clear that under the
assumption of part (i) of Theorem 2.2, one has
\[
\int G_{\alpha V}(x, t; y, s) dy \le C <\infty.
\]This is one of the main assumptions used by Sturm [St] in the
time independent case (Theorem 4.12).
  If $\alpha=1$, then the conclusion of Theorem 2.2 may not
hold even for time independent potentials. See [St]. Also note that
this theorem provides a nearly necessary and sufficient condition on
certain sign changing potential such that the corresponding heat
kernel has Gaussian upper and lower bound. The only "gap" in the
condition is the difference in the parameters of the kernels
$g_{1/4}$ and $g_b$. It is well know and easy to check if $|\nabla
f| \in L^{p, q}$ with $\frac{n}{p} + \frac{2}{q}<1$ and $f=0$
outside a compact set, then $g_b \star |\nabla f|^2$ is a bounded
function for all $b>0$.

\begin{corollary} Let $V \in L^2(\R^n \times (0, \infty))$.

(a). Suppose
\begin{equation}
\label{domi} \int^T_0 \int  V \phi^2 \le \int^T_0 \int |\nabla
\phi|^2 + b \int^T_0 \int \phi^2.
\end{equation}
for all smooth,
compactly supported function $\phi$ in $\R^n \times (0, T)$ and some
$b>0$ and $T>0$. Then (\ref{1.1}) has a positive solution when $u_0
\ge 0$ and moreover
\[
V = \Delta f - |\nabla f|^2 -\partial_t f.
\]

(b). Suppose $V = \Delta f - |\nabla f|^2$ then $V$ is form bounded,
i.e. it satisfies (\ref{domi}).

(c). Let $V$ be a $L^2$ potential permitting positive $L^2$
solutions for (1.1), then $V$ can be written as the sum of one form
bounded potentials and the time derivative of a function almost
bounded from above by a constant.
\end{corollary}
\medskip

In the next corollary, we consider only time-independent,
nonnegative potentials. Here the definition of $V = \Delta f -
|\nabla f|^2$ is slightly different from that of Definition (1.1)
since we do not have to worry about time derivatives. One
interesting consequence is that these class of potentials is {\it
exactly} the usual form boundedness potentials.

\begin{corollary} Suppose
$0 \le V \in L^1(\R^n)$.
 Then the following statements are equivalent.

(1). For some $f \in L^1_{loc}(\R^n)$ and a constant $b>0$,
\[
V = \Delta f - |\nabla f|^2 + b.
\]This mean there exist $V_j \in L^\infty$ such that $V_j \to V$ in
$L^2(\R^1)$ and $V_j = \Delta f_j - |\nabla f_j|^2 +b$ for some $f_j
\in W^{2, 2}(\R^n)$, \ $j=1, 2, ...$.

(2).
\[
\int  V \phi^2 dx \le \int |\nabla \phi|^2 dx + b  \int \phi^2
dx.
\]for all smooth, compactly supported function $\phi$ in $\R^n$ and some
$b \ge 0$.
\end{corollary}
\medskip

{\it Remark 2.3.} Condition (2) in the above corollary just means
that the bottom of the spectrum for the operator $-\Delta -V$ is
finite. This condition is the same as those given in [CM] and [GZ].

It is a fact that most people feel more familiar with the case when
the potential $V$ is written as $L^{p, q}$ functions. Also there may
be some inconvenience about the presence of the nonlinear term
$|\nabla f|^2$ in the potential in Theorem 2.2. Therefore in our
next theorem, we will use only $L^{p, q}$ conditions on $f$ without
nonlinear terms.

\begin{theorem}
Suppose

(a) $f \in L^\infty(\R^n \times \R^+)$;

(b) $f=0$ outside a cylinder $B(0, R_0) \times [0, T_0]$, \ $R_0,
T_0>0$;

(c) $|\nabla f|\in L^{p, q}(\R^n \times \R^+)$ with $\frac{n}{p} +
\frac{2}{q}<1$;

(d) $V = \Delta f - \partial_t f \in L^2(\R^n \times \R^+)$.

Then there exists a constant $A_0$ depending only on $n, p, q$ such
that the following statements hold, provided that
\[
\Vert \  |\nabla f| \  \Vert_{L^{p, q}(\R^n \times \R^+)} < A_0.
\]

(i) The kernel $G_V$ has Gaussian upper and lower bound in all space
time.

(ii) The kernel $G_{\partial_t f}$ has Gaussian upper and lower
bound in all space time.
\end{theorem}
\medskip

{\it Remark 2.4.}  If $V$ is independent of time, then Theorem 2.3.
reduces to the well known classical fact:

 if a potential $V$ is the derivative a of a small $L^{n +\epsilon}$
 function, then $G_V$ has Gaussian upper and lower bound. (See
 [LS]) e.g.

 In the time dependent case our result is genuinely new due to the
 presence of the term $\partial_t f$. Let us mention that some
 smallness condition on the potential is needed for the existence of
 Gaussian bounds for $G_V$. This is the case even for time
 independent, smooth potentials due to the possible presence of
 ground state.

\medskip
\subsection{Preliminaries}

In order to prove the theorems we need to prove a proposition
concerning the existence, uniqueness and maximum principle for
solutions of (\ref{1.1}) under the assumptions that $V$ is bounded
from above by a constant. The result is standard if one assumes that
the gradient of solutions are $L^2$. However we only assume that
solutions are $L^2$. Therefore a little extra work is needed.

\medskip

\begin{proposition}
Suppose that $V \in L^2(\R^n \times \R^+)$ and that $V \le b$ for a
positive constant $b$. Then the following conclusions hold.

(i). The only $L^2$ solution to the problem
\[
\begin{cases}
\Delta u + V u - \partial_t u = 0, \quad \R^n \times (0, T), \
T>0, \\
u(x, 0) = 0,
\end{cases}
\]
is zero.

(ii). Let $u$ be a solution to the problem in (i) such that
$u(\cdot, t) \in L^1(\R^n)$, \  $u \in L^1(\R^n \times (0, T))$ and
$V u \in L^1(\R^n \times (0, T))$.  Then $u$ is identically zero.

(iii). Under the same assumptions as in (i), the following problem
has a unique $L^2$ nonnegative solution.
\[
\begin{cases}
\Delta u + V u - \partial_t u = 0, \quad \R^n \times (0, T), \
T>0, \\
u(\cdot, 0) = u_0(\cdot)\ge 0, \quad u_0 \in L^2(\R^n).
\end{cases}
\]

(iv). Under the same assumptions as in (i), let $u$ be a $L^2$
solution to the following problem
\[
\begin{cases}
\Delta u + V u - \partial_t u = f, \quad \R^n \times (0, T), \
T>0, \\
u(\cdot, 0) = 0.
\end{cases}
\]Here $f \le 0$ and $f \in L^1(\R^n \times (0, T))$. Then $u \ge
0$ in $\R^n \times (0, T)$.
\end{proposition}

Proof of (i). Let $u$ be a $L^2$ solution to the problem in (i).
Choose a standard mollifier $\rho$ and define, for $j = 1, 2, ...$,
\[
u_j(x, t) = j^n \int \rho( j (x-y)) u(y, t) dy \equiv \int \rho_j
(x-y) u(y, t) dy.
\]Then $\nabla u_j$ and $\Delta u_j$ exist in the classical sense.
From the equation on $u$, it holds
\[
\Delta u_j + \int \rho_j (x-y) V(y, t) u(y, t) dy - \partial_t u_j =
0,
\]where $\partial_t u_j$ is understood in the weak sense.

Given $\epsilon>0$, we define
\[
h_j = \sqrt{ u^2_j + \epsilon}.
\]Let $0 \le \phi \in C^\infty_0(\R^n)$. Then direct calculation shows
\[
\begin{aligned}
&\int \phi h_j(x, s) |^t_0 dx = \\
&\int^t_0\int \frac{u_j \Delta u_j}{\sqrt{ u^2_j + \epsilon}}(x, t)
\phi(x) dxds + \int^t_0\int \phi(x) \frac{u_j(x, s)}{\sqrt{
u^2_j + \epsilon}} \ \int \rho_j (x-y) V(y, t) u(y, t) dy dxdt\\
&\equiv T_1 + T_2.
\end{aligned}
\]Using integration by parts, we deduce
\[
\begin{aligned}
T_1 &= - \int^t_0\int \frac{|\nabla u_j |^2}{\sqrt{ u^2_j +
\epsilon}} \phi(x, s) dxds + \int^t_0\int \frac{u^2_j}{u^2_j +
\epsilon} \frac{|\nabla u_j |^2}{\sqrt{ u^2_j + \epsilon}} \phi(x,
s) dxds\\
&\qquad - \int^t_0\int u_j \frac{\nabla u_j \nabla \phi}{\sqrt{
u^2_j + \epsilon}}  dxds
\end{aligned}
\]Since the sum of the first two terms on the righthand side of
the above inequality is non-positive, we have
\[
T_1 \le \int^t_0\int u_j \frac{|\nabla u_j| \  |\nabla \phi|}{\sqrt{
u^2_j + \epsilon}}  dxds.
\]Taking $\epsilon$ to zero, we obtain
\[
\begin{aligned}
\int \phi |u_j(x, t)|  dx & \\
&\le \int^t_0\int |\nabla u_j|  |\nabla \phi| dxds + \int^t_0\int
\phi(x) | \int \rho_j V^+ u(y, s) dy | dxds \\
&\qquad - \int^t_0\int \phi(x) \frac{u_j}{|u_j|}(x, s) \int
\rho_j(x-y) V^- u(y, s) dy  dxds.
\end{aligned}
\]Here and later we set $\frac{u_j}{|u_j|}(x, s)=0$ if $u_j(x, s)
= 0$.

Next, since $u, V \in L^2(\R^n \times (0, T))$, one has $u V \in
L^1(\R^n \times (0, T))$. From the equation
\[
\Delta u + V u - \partial_t u = 0
\]one deduces
\[
u(x, t) = \int^t_0\int G_0(x, t; y, s) (V u)(y, s) dyds.
\]Here, as always, $G_0$ is the fundamental solution of the
free heat equation. Hence $u(\cdot, t) \in L^1(\R^n)$ and $u \in
L^1(\R^n \times (0, T))$. Therefore, for any fixed $j$, there holds
\[
|\nabla u_j | \in L^1(\R^n \times (0, T)).
\]Now, for each $R>0$,  we choose $\phi$ so that $\phi = 1$ in
$B(0, R)$, $\phi=0$ in $B(0, R+1)^c$ and $|\nabla \phi| \le 2$.
Observing
\[
\int^t_0 \int |\nabla u_j | \ |\nabla \phi| dxds \to 0, \qquad R \to
\infty,
\]we deduce, by letting $R \to \infty$,
\begin{equation}
\label{uj}
\int  |u_j(x, t)|  dx  \\
\le  \int^t_0\int
 | \int \rho_j V^+ u(y, s) dy | dxds \\
- \int^t_0\int  \frac{u_j}{|u_j|}(x, s) \int \rho_j(x-y) V^- u(y, s)
dy  dxds.
\end{equation}By the fact that $V^-u, V^+ u \in L^1(\R^n
\times (0, T))$, we know that
\[
\int \rho_j(\cdot-y) V^- u(y, \cdot) dy \to V^- u(\cdot, \cdot),
\]
\[
\int \rho_j(\cdot-y) V^+ u(y, \cdot) dy \to V^+ u(\cdot, \cdot),
\]in $L^1(\R^n \times (0, T))$.
Since $\frac{u_j}{|u_j|}(x, s)$ is bounded and converges to
$\frac{u}{|u|}(x, s)$ a.e. in the support of $u$, we have
\[
\begin{aligned}
& \bigg{|} \int^t_0\int \big{[} \frac{u_j}{|u_j|}(x, s) \int
\rho_j(x-y) V^- u(y, s) dy   - \frac{V^- u^2}{|u|} \big{]} dxds
\bigg{|}\\
&\le \bigg{|} \int^t_0\int \big{(}  \frac{u_j}{|u_j|} -
\frac{u}{|u|} \big{)} V^-  u dxds \bigg{|} + \bigg{|} \int^t_0\int
\frac{u_j}{|u_j|} \big{(}  \rho_j \star V^-u - V^- u  \big{)}dxds
\bigg{|} \to 0.
\end{aligned}
\]

Substituting this to (\ref{uj}) we deduce, by taking $j \to \infty$,
\[
\int  |u(x, t)|  dx  \\
\le  \int^t_0\int V^+ |u(x, s)|  dxds  - \int^t_0\int \frac{V^-
u^2}{|u|}(x, s) dxds.
\]Therefore
\[
\int  |u(x, t)|  dx  \\
\le  \int^t_0\int  |u(x, s)|  dxds \ \Vert V^+ \Vert_{\infty}.
\]By Grownwall's inequality $u(x, t) =0$ a.e. This proves part
(i).
\medskip

{\it Proof of (ii).}

Notice that the only place we have used the $L^2$ boundedness of $u$
is to ensure that $V u \in L^1(\R^n \times (0, T))$. But this a part
of the assumptions in (ii). Therefore (ii) is also proven.

\medskip

{\it Proof of (iii).}  The uniqueness is an immediate consequence of
part (i). So we only need to prove existence. This follows from a
standard limiting process. For completeness we sketch the proof.

Given $k=1, 2, ..$ let $V_k$ be the truncated potential
\[
V_k = \sup \{ V(x, t), - k \}.
\]
Since $V_k$ is a bounded function there exists a unique, nonnegative
solution $u_k$ to the following problem.
\[
\begin{cases}
 \Delta u_k + V_k u_k - \partial_t u_k=0, \quad \text{in} \quad {\bf R}^n \times
(0, \infty),  \\
u(x, 0) = u_0(x) \ge 0, \quad x \in {\bf R}^n, \quad u_0 \in
L^2(\R^n).
\end{cases}
\]
By the standard maximum principle, $\{ u_k \}$ is a nonincreasing
sequence and
\[
\frac{1}{2} \int u^2_k |^T_0 dx + \int^T_0\int |\nabla u_k |^2dxdt =
\int^T_0\int v_k u^2_k dxdt \le \int^T_0\int V^+_k u^2_k dxdt\le
\Vert V^+ \Vert_\infty \int^T_0\int  u^2_k dxdt.
\]By Grownwall's lemma, we have
\[
\int^T_0\int |\nabla u_k |^2 dxdt + \int u^2_k(x, T) dx \le \int
u^2_0(x) dx + \Vert v^+ \Vert_\infty \int u^2_0(x) dx \ e^{2 \Vert
v^+ \Vert_\infty T}.
\]It follows that $u_k$ converges pointwise to a function $u$ which
also satisfies the above inequality. Let $\phi$ be a test function
with compact support. Then
\[
\int ( u_k  \phi)  |^{T}_0 dx -\int^{T}_0 \int u_k \phi_t dxdt -
\int^{T}_0 \int  u_k  \Delta \phi dxdt - \int^T_0 \int V_k u_k \phi
dxdt =0
\]Since $\{ u_k \}$ is a monotone sequence and also since $|V_k u_k|
\le |V| u_1 \in L^1(\R^n \times (0, T))$, the dominated convergence
theorem implies that
\[
\int ( u  \phi)  |^{T}_0 dx -\int^{T}_0 \int u \phi_t dxdt -
\int^{T}_0 \int  u  \Delta \phi dxdt - \int^T_0 \int V u \phi dxdt
=0.
\]This shows that $u$ is a nonnegative solution. It is clear that
$u$ is not identically zero since $u_0$ is not. This proves part
(iii) of the proposition.
\medskip

{\it Proof of Part (iv).}

Let $V_k$ be a truncated potential as in part (iii). Since $V_k$ is
bounded, the standard maximum principle shows that there exists a
unique, nonnegative solution to the following problem.
\[
\begin{cases}
\Delta u_k + V_k u_k - \partial_t u_k = f \le 0, \quad \R^n \times
(0, T), \
T>0, \\
u_k(\cdot, 0) = 0.
\end{cases}
\]Moreover $\{ u_k \}$  forms a decreasing sequence. Since $V_k$ is a
bounded function, the standard parabolic theory shows that
\[
\int u_k(x, t) dx = \int^t_0 \int V^+_k u_k dxds - \int^t_0 \int
V^-_k u_k dxds  + \int^t_0 \int f dx ds.
\]Therefore
\[
\int u_k(x, t) dx \le \Vert V^+ \Vert_{\infty} \int^t_0 \int u_k
dxds + \int^t_0 \int f dx ds.
\]This implies
\[
\int u_k(x, t) dx +  \int^t_0 \int u_k dxds \le C(t, \Vert V^+
\Vert_{\infty}, \Vert f \Vert_1).
\]It follows that
\[
\int u_k(x, t) dx +  \int^t_0 \int u_k dxds + \int^t_0 \int V^-_k
u_k dxds \le C(t, \Vert V^+ \Vert_{\infty}, \Vert f \Vert_1).
\]Let $w$ be the pointwise limit of the decreasing sequence $\{ u_k
\}$. Then we have
\[
\int w(x, t) dx +  \int^t_0 \int w dxds \le C(t, \Vert V^+
\Vert_{\infty}, \Vert f \Vert_1).
\]It is straight forward to check that $w$ is a nonnegative solution
to the problem
\[
\begin{cases}
\Delta w + V w - \partial_t w = f \le 0, \quad \R^n \times (0, T), \
T>0, \\
w(\cdot, 0) = 0.
\end{cases}
\]Hence
\[
\begin{cases}
\Delta (w-u) + V (w-u) - \partial_t (w-u) =  0, \quad \R^n \times
(0, T), \
T>0, \\
(w-u)(\cdot, 0) = 0.
\end{cases}
\]Recall that $u$ is assumed to be a $L^2$ solution and that $V \in
L^2$. We have that $V u \in L^1$ and consequently $u(\cdot, t) \in
L^1(\R^n)$ and $u \in L^1(\R^n \times (0, T))$. Now by Part (ii) of
the proposition, we deduce $w=u$ since $w$ is also $L^1$. Hence $u
\ge 0$. This finishes the proof of the proposition. \qed

\medskip

\subsection{Proofs of Theorems}

\medskip
\[
\]
{\it Proof of Theorem 2.1 (i)}.

For $j=1, 2...$, let
\[
V_j= \min \{V(x, t), j \}.
\]Since $V_j$ is bounded from above, Proposition 2.1 shows that
there exists a unique solution $u_j$ to the following problem.
\[
\begin{cases}
\Delta u_j + V_j u_j -\partial_t u_j = 0, \qquad (x, t) \in \R^n \times (0, \infty)\\
u_j(x, 0) = u_0(x), \qquad x \in \R^n.
\end{cases}
\]Notice that $u_j - u_{j-1}$ is a solution to the problem
\[
\begin{cases}
\Delta (u_j-u_{j-1}) + V_j (u_j-u_{j-1}) -\partial_t (u_j-u_{j-1}) =
(V_{j-1}-V_j) u_{j-1}, \qquad (x, t) \in \R^n \times (0, \infty)\\
(u_j-u_{j-1})(x, 0) = 0, \qquad x \in \R^n.
\end{cases}
\]Notice that
\[
(V_{j-1}-V_j) u_{j-1} \le 0, \qquad  (V_{j-1}-V_j) u_{j-1} \in
L^1(\R^n \times (0, T)), \quad T>0.
\]We can then apply Proposition 2.1 (iv) to conclude that
\[
u_j \ge u_{j-1}.
\]Moreover
\[
\begin{cases}
\Delta (u-u_j) + V_j (u-u_j) -\partial_t (u-u_j) =
(V_j-V) u, \qquad (x, t) \in \R^n \times (0, \infty)\\
(u-u_j)(x, 0) = 0, \qquad x \in \R^n,
\end{cases}
\]with
\[
(V_j-V) u \le 0, \qquad  (V_j-V) u \in L^1(\R^n \times (0, T)),
\quad T>0.
\]By Proposition 2.1 (iv) again we know that
\[
u \ge u_j.
\]Therefore $\{ u_j \}$ is a non-decreasing sequence of nonnegative functions
bounded from above by a $L^2$ function. Let $w$ be the pointwise
limit of $u_j$. The $w$ is $L^2$ and $|V_j u_j| \le |V u| \in
L^1(\R^n \times (0, T))$, $T>0$. By the dominated convergence
theorem, it is straight forward to check that $w$ is a nonnegative
$L^2$ solution to the equation
\[
\begin{cases}
\Delta w + V w -\partial_t w = 0, \qquad (x, t) \in \R^n \times (0, \infty)\\
w(x, 0) = u_0(x), \qquad x \in \R^n.
\end{cases}
\]

Fixing $j$, for any $k=1, 2, ...$, Let
\[
V_{jk} = \max \{ V_j(x, t), -k \}.
\]Since $V_{jk}$ is bounded, the following problem has a unique
$L^2$ solution.
\[
\begin{cases}
\Delta u_{jk} + V_{jk} u_{jk} -\partial_t u_{jk} = 0,
\qquad (x, t) \in \R^n \times (0, \infty)\\
u_j(x, 0) = u_0(x), \qquad x \in \R^n.
\end{cases}
\]Due to the fact that $\{ V_{jk} \}$ is a decreasing sequence of
$k$, the maximum principle shows that $\{ u_{jk} \}$ is also a
decreasing sequence of $k$. Since $V_{jk} u_{jk} \in L^2(\R^n \times
(0, T))$, $T>0$, the parabolic version of the Calderon-Zygmond
theory shows
\[
\Delta u_{jk}, \quad \partial_t u_{jk} \in L^2(\R^n \times (0, T)),
T>0.
\]

Since
\[
0 \le u_{jk} - u_j \le u_{j1}- u_j, \ k =1, 2, 3, ...,
\]
\[
0 \le w - u_j \le w-u_1,
\]we can apply the dominated convergence theorem to conclude that
\[
\lim_{k \to \infty} \int^T_0 \int ( u_{jk} -u_j)^2 dxdt = 0, \quad
\lim_{j \to \infty} \int^T_0 \int (w -u_j)^2 dxdt=0.
\]Therefore we can extract a subsequence $\{ u_{jk_j} \}$ such that
\[
\lim_{j \to \infty} \int^T_0 \int (w -u_{jk_j})^2 dxdt=0.
\]Hence there exists a subsequence, still called $\{ u_{jk_j} \}$ such that
\[
 u_{jk_j} \to w \qquad a.e.
\]

Recall that $u_0 \ge 0$,  $u_0 \neq 0$ and $V_{jk_j}$ is bounded. It
is clear that $u_{jk_j} >0$ when $t>0$. Now we define
\[
f_j = - \ln u_{jk_j}, \qquad f = - \log w.
\]
Then
\[
V_{jk_j} = \Delta f_j - |\nabla f_j|^2 - \partial_t f_j.
\]Clearly $f_j \to f$ a.e. and $V_{jk_j} \to V$ in $L^2$ as $j \to
\infty$. By Definition 2.1, this means
\[
V = \Delta f - |\nabla f|^2 - \partial_t f.
\]It is clear that $e^{-f}b = w$ is $L^2$ by construction.
\medskip

{\it Proof of Theorem 2.1 (ii).}

By assumption, there exist sequences of functions $\{ V_j \}$ and
$\{ f_j \}$ such that
\[
\Vert V_j \Vert_{L^2} \le C, \quad |V_j| \le |V_{j+1}|, \quad f_j
\in L^\infty, \quad \Vert V_j -V \Vert_{L^2} \to \infty,
\]
\[
V_j = \Delta f_j - |\nabla f_j|^2 -\partial_t f_j.
\]Then for $u_j \equiv e^{-f_j}$, we have
\[
\Delta u_j + V_j u_j -\partial_t u_j = 0.
\]We will show that $\Vert u_j \Vert_{L^2}$ is uniformly bounded.
To this end, we observe that
\[
\Delta (u_j - u_{j+1}) + V_j (u_j - u_{j+1}) - \partial_t (u_j -
u_{j+1}) = - (V_j -V_{j+1}) u_{j+1} \ge 0.
\]Recall from Definition 2.1 that $f_j(x, 0)$ is independent of $j$. Hence
$u_j(x, 0) = u_{j+1}(x, 0)$. Therefore $0 \le u_j \le u_{j+1}$. By
the assumption that $f_j \to f$ a.e., we know that $u_j = e^{- f_j}
\to e^{-f}$ a.e. Note that $ e^{-f} \in L^2(\R^n \times (0,
\infty))$. Hence $\Vert u_j \Vert_{L^2}$ is uniformly bounded.

By weak compactness in $L^2$, there exists a subsequence, still
called $\{ u_j \}$ such that $u_j$ converges weakly to a $L^2$
function which we will call $u$.  Observe that, for any compactly
supported test function $\phi$, there holds
\[
\begin{aligned}
&\Vert u_j V_j \phi - u V \phi \Vert_{L^1} \le \Vert u_j (V_j - V)
\phi \Vert_{L^1} + \Vert (u_j-u) V  \phi \Vert_{L^1}\\
&\le \Vert u_j \Vert_{L^2} \Vert V_j -  V \Vert_{L^2} \ \Vert \phi
\Vert_{L^\infty} + \Vert (u_j-u) V  \phi \Vert_{L^1}.
\end{aligned}
\]Hence
\[
\Vert u_j V_j \phi - u V \phi \Vert_{L^1} \to 0
\]when $j \to \infty$. From here it is easy to check that $u$ is a
nonnegative solution to (\ref{1.1}) with $u_0 = e^{-f_j(x, 0)}$ as
the initial value. Note the $u_0$ is independent of $j$. If $u_0 \in
L^2(\R^n)$ then we are done. Otherwise, we can selection a $L^2$
function dominated by $u_0$ to serve as the initial value. \qed
\medskip

Next we will provide a

\bigskip
{\it Proof of Theorem 2.2 (i)}.
\medskip

We will use an idea based on an argument  in [St] where the heat
equation with some singular, time independent potentials are
studied.

By virtue of Proposition 2.1, the fundamental solution $G_V$ is
defined as the limit of fundamental solutions of the equation in
(\ref{1.1}) where $V$ is replaced by nonsingular potentials.
Therefore we can and will assume that $V$ is smooth in this
subsection. The constants involved will be independent of the
smoothness.

Since, by assumption
\beq
 \label{Vf}
  V = \Delta f - \alpha |\nabla f|^2 -
\partial_t f,
\enq
one has
\[
\alpha V = \Delta (\alpha f) -  |\nabla (\alpha f)|^2 - \partial_t
(\alpha f).
\]Writing $F = e^{- \alpha f}$, it is easy to show that
\beq \label{alphaV}
 \Delta F + \alpha V F - \partial_t F = 0.
\enq Let us denote the fundamental solution of the equation in
(\ref{alphaV}) by $G_{\alpha V}$. Since $f$ is bounded, we know that
$F$ is bounded between two positive constants. Therefore it is clear
that
\beq \label{IG_a}
 0< \frac{\inf F}{\sup F} \le \int G_{\alpha
V}(x, t; y, s) dy \le \frac{\sup F}{\inf F}, \enq for all $x \in
\R^n$ and $t>s$. Here $\inf F$ and $\sup F$ are taken over the whole
domain of $F$.

By Feynman-Kac formula and H\"older's inequality, for a given $\phi
\in C^\infty_0(\R^n)$, there holds
\[
\big{|} \int  G_V (x, t; y, s) \phi(y) dy \big{|} \le \big{[} \int
G_{\alpha V} (x, t; y, s)  dy \big{]}^{1/\alpha} \ \big{[} \int G_0
(x, t; y, s) |\phi(y)|^{\alpha/(\alpha-1)} dy
\big{]}^{(\alpha-1)/\alpha}
\]By (\ref{IG_a}), we deduce
\[
\big{|} \int  G_V (x, t; y, s) \phi(y) dy \big{|} \le \frac{c
s^{1/\alpha}_0}{(t-s)^{(\alpha-1)n/(2\alpha)}} \Vert \phi
\Vert_{\alpha/(\alpha-1)},
\]where $s_0 = \frac{\sup F}{inf F}$. The norm on $\phi$ means the
$L^{\alpha/(\alpha-1)}(\R^n)$ norm.
Hence
\begin{equation}
\label{Gpq}
 \Vert G_V (\cdot, t; \cdot, s) \Vert_{\alpha/(\alpha-1),
\infty} \le \frac{c s^{1/\alpha}_0}{(t-s)^{(\alpha-1)n/(2\alpha)}}.
\end{equation}
Here and later the norm $\Vert \cdot \Vert_{p, q}$ stands for the
operator norm from $L^p(\R^n)$ to $L^q(\R^n)$ for $p, q$ between $1$
and $\infty$.

Without loss of generality we assume that $\alpha/(\alpha-1)$ is an
integer. This is so because otherwise we can choose one $\alpha_1
\in (1, \alpha)$ such that $\alpha_1/(\alpha_1-1)$ is an integer.
Then interpolating between $G_{\alpha V}$ and $G_0$ by Feynman-Kac
formula again, we know that
\[
\int G_{\alpha_1 V}(x, t; y, s) dy \le C(F, \alpha, \alpha_1).
\]Then we can just work with $G_{\alpha_1 V}$ instead of $G_{\alpha
V}$ in the above.

Using the reproducing property of $G_V$ we deduce
 \beq \label{G1}
\Vert G_V (\cdot, t; \cdot, s) \Vert_{1, \infty} \le \Pi^m_{j=1}
\Vert G_V (\cdot, s + t_j; \ \cdot, t_{j-1}) \Vert_{p_j, q_j},
\enq
where
\[
m = \alpha/(\alpha-1), \quad p_j=m/(m-l+1), \quad q_j=m/(m-l), \quad
t_j= s + ((t-s)j/m).
\]For each $j$ between $1$ and $m$, we apply the Riesz-Thorin
interpolation theorem to deduce
\[
\Vert G_V (\cdot, t_j; \ \cdot, t_{j-1}) \Vert_{p_j, q_j} \le \Vert
G_V (\cdot, s + t_j; \ \cdot,  t_{j-1}) \Vert^{1-\lambda_j}_{1,
m/(m-1)}\ \Vert G_V (\cdot, t_j; \ \cdot,  t_{j-1})
\Vert^{\lambda_j}_{m, \infty}.
\]Here the parameters are determined by the following relations
\[
\frac{1}{p_j} = \frac{1-\lambda_j}{1} + \frac{\lambda_j}{m}, \qquad
\frac{1}{q_j} = \frac{1-\lambda_j}{m/(m-1)} +
\frac{\lambda_j}{\infty},
\]
\[
\lambda_j = \frac{j-1}{m-1}.
\]It follows that
\begin{equation}
\label{Gpj}
 \Vert G_V (\cdot, t_j; \ \cdot, t_{j-1}) \Vert_{p_j,
q_j} \le \Vert G_V (\cdot, t_j; \ \cdot,  t_{j-1}) \Vert_{m,
\infty}.
\end{equation}
Substituting (\ref{Gpq}) to (\ref{Gpj}), we deduce, after noticing
that $t_j-t_{j-1} = (t-s)/m$,
\[
\Vert G_V (\cdot, t_j; \ \cdot, t_{j-1}) \Vert_{p_j, q_j} \le
\frac{c s^{1/\alpha}_0
m^{(\alpha-1)n/(2\alpha)}}{(t-s)^{(\alpha-1)n/(2\alpha)}}.
\]This and (\ref{G1}) imply that
\[
\Vert G_V (\cdot, t; \cdot, s) \Vert_{1, \infty} \le \frac{c
s^{m/\alpha}_0 m^{n/2}}{(t-s)^{n/2}}.
\]Here we just used the relation $m = \alpha/(\alpha-1)$. This yields
the on-diagonal upper bound \beq \label{odub}
 G_V(x, t; y, s)  \le \frac{c s^{1/(\alpha-1)}_0
(\alpha/(\alpha-1))^{n/2}}{(t-s)^{n/2}}. \enq In order to obtain the
full Gaussian bound, we observe that, for any $p>1$, the Feynman-Kac
formula implies
\beq \label{FKp} G_V(x, t; y, s) \le \big{[}
G_{pV}(x, t; y, s) \big{]}^{1/p} \ \big{[} G_0(x, t; y, s)
\big{]}^{(p-1)/p}.
\enq
Notice also
\[
p V = p (\Delta f - \alpha |\nabla f|^2 - \partial_t f) = \Delta (p
f) - \frac{\alpha}{p} |\nabla (p f)|^2 - \partial_t(p f).
\]Taking $p = (1+\alpha)/2$, then $\frac{\alpha}{p}>1$. Therefore
$p V$ also satisfies the condition of Theorem 2.2 (i). Hence, the
on-diagonal bound (\ref{odub}) holds for $G_{pV}$. i.e., there
exists a constant $C(\alpha, e^{\sup f-\inf f})$ such that
\[
 G_{p V}(x, t; y, s)  \le \frac{C(\alpha, e^{\sup f-\inf f})}{(t-s)^{n/2}}.
\]Substituting this to the inequality (\ref{FKp}), we obtain
the desired Gaussian upper bound for $G_V$.

\medskip

{\it Proof of (ii).}

 In this part we prove the Gaussian lower bound.
We will follow Nash's original idea. The novelty is a way of
handling the potential term even if it is very singular. The main
idea is to exploit the structure of the potential when it is written
as a combination of derivatives.

Since the setting of our problem is invariant under the scaling, for
$r>0$,
\[
V_r(x, t) = r^2 V(r x, r^2 t), \qquad f_r(x, t) = f(r x, r^2 t),
\qquad u_r(x, t) = r^2 u(r x, r^2 t),
\]we can just prove the lower bound for $t=1$ and $s=0$. We divide
the proof into three steps.

{\it Step 1.} Fixing $x \in \R^n$, let us set
\[
u(y, s) = G_V(y, s; x, 0),
\]
\[
H(s) = \int e^{- \pi |y|^2} \ln u(y, s) dy.
\]Differentiating $H(s)$, one obtains
\[
\begin{aligned}
H'(s) &= \int e^{- \pi |y|^2} \frac{\partial_s u(y, s)}{u(y, s)}
dy\\
&= - \int \nabla \big{(} \frac{e^{- \pi |y|^2}}{u} \big{)} \nabla u
dy + \int e^{- \pi |y|^2} V(y, s) dy.
\end{aligned}
\]Estimating the first term on the righthand side of the above
inequality as in [FS], section 2, one arrives at
\begin{equation}
\label{H'(s)}
 H'(s) \ge  - C + \frac{1}{2} \int e^{- \pi |y|^2}
|\nabla \ln u(y, s)|^2 dy +
  \int e^{- \pi
|y|^2} V(y, s) dy.
\end{equation}
Here $C$ is a positive
constant. Since,
\[
V = \Delta f - \alpha |\nabla f|^2 - \partial_t f,
\]we know that
\[
\begin{aligned}
f(x, t)& = \int G_0(x, t; y, 0) f(y, 0)\\
&\qquad - \int^t_0 \int G_0(x, t; y, s) V(y, s) dyds - \alpha
\int^t_0 \int G_0(x, t; y, s) |\nabla f|^2 dyds.
\end{aligned}
\]By our assumption
\[
\int^t_0 \int G_0(x, t; y, s) |\nabla f|^2 dyds<\infty.
\]Hence the boundedness of $f$ implies that
\[
m(x, t) \equiv - \int^t_0 \int G_0(x, t; y, s) V(y, s) dyds \in
L^\infty.
\]Moreover
\beq
\label{MV}
 \Delta m - \partial_t m = V. \enq Therefore
\[
\begin{aligned}
\int e^{- \pi |y|^2} &V(y, s) dy
=\int e^{- \pi |y|^2} [\Delta m - \partial_s m] (y, s) dy\\
&=\int [\Delta e^{- \pi |y|^2}]  m dy  - \partial_s \int e^{- \pi
|y|^2} m(y, s) dy\\
&\ge - C - \partial_s \int e^{- \pi |y|^2} m(y, s) dy.
\end{aligned}
\]Here we have used the boundedness of $m$. Substituting the above
to the righthand side of (\ref{H'(s)}), we obtain
\begin{equation}
\label{H'M'} H'(s) \ge  - C + \frac{1}{2} \int e^{- \pi |y|^2}
|\nabla \ln u(y, s)|^2 dy - M'(s),
\end{equation}where
\begin{equation}
\label{M(s)}
 M(s) = \int e^{- \pi |y|^2} m(y, s) dy.
\end{equation}

\medskip

{\it Step 2.} By Poincar\'e's inequality with $e^{- \pi |y|^2}$ as
weight, we deduce, for some $B>0$,
\[
H'(s) \ge -C + B \int e^{- \pi |y|^2} [ \ln u(y, s) - H(s) ]^2 dy -
M'(s).
\]Next, observe that $(\ln u - H(s))^2/u$ is  non-increasing as a
function of $u$ when $u$ is between $e^{2+H(s)}$ and $\infty$. Also
from the Gaussian upper bound,
\[
\sup_{1/2 \le s \le 1} u(y, s) \le K < \infty.
\]Therefore
\beq
\label{H'}
 H'(s) \ge -C + C B  K^{-1} (\ln K - H(s))^2 \int_{u(y, s) \ge
\exp(2+H(s))} e^{- \pi |y|^2} u(y, s) dy - M'(s).
\enq

Using the Gaussian upper bound again, we know that $H(s) \le C$ for
some $C>0$ and that
\begin{equation}
\label{UH}
\begin{aligned}
 \int_{u(y, s) \ge exp(2+H(s))}& e^{- \pi |y|^2} u(y, s) dy
\ge \int e^{- \pi |y|^2} u(y, s) dy - c e^{2 + H(s)}\\
&\ge e^{- \pi r^2} \int_{|y|<r} u(y, s) dy - c e^{2+H(s)}\\
&=e^{- \pi r^2} \big{[} \int u(y, s) dy  - \int_{|y|<r} u(y, s) dy
\big{]}- c e^{2+H(s)}.
\end{aligned}
\end{equation}

We aim to find a lower bound for the righthand side of (\ref{UH}).
By (\ref{MV}),
\[
V = \Delta m - \partial_t m \ge  \Delta m - |\nabla m |^2 -
\partial_t m \equiv V_1.
\]Write $h = e^{-m}$. Then
\[
\Delta h + V_ 1 h - \partial_t u = 0.
\]Since $m \in L^\infty$, we know that $h$ is bounded between two
positive constants. Observe that
\[
h(x, t) = \int G_{V_1}(x, t; y, s) h(y, s) dy.
\]Hence
\[
0<c_1 < \int G_{V_1}(x, t; y, s) dy \le c_2.
\]By the maximum principle, we have
\beq \label{ILB} \int G_V(x, t; y, s) dy \ge \int G_{V_1}(x, t; y,
s) dy \ge c_1>0. \enq

Recall that $u(y, s) = G_V(y, s; x, 0)$. Substituting (\ref{ILB}) to
(\ref{UH}) and applying the Gaussian upper bound on $u(y, s)$, we
deduce
\begin{equation}
\label{CUH}
 \int_{u(y, s) \ge exp(2+H(s))} e^{- \pi |y|^2} u(y, s) dy
\ge  c e^{- \pi r^2} c_1- c e^{2+H(s)},
\end{equation}when $r$ is sufficiently large. Substituting (\ref{CUH}) to
(\ref{H'}), we arrive at
\beq
\label{H'2}
 H'(s) \ge -C + C B  K^{-1} (\ln K - H(s))^2   e^{- \pi r^2} c_1- c
 e^{2+H(s)} - M'(s).
\enq

We claim that $H(1) \ge -c_0$ for some sufficiently large $c_0>0$.
Suppose otherwise, i.e. $H(1) < - c_0$.  From (\ref{H'2}), for some
$C>0$,
\[
H'(s) \ge -C  - M'(s).
\]Hence
\[
H(1)-H(s) \ge -C (1-s) - (M(1)-M(s)).
\]Therefore
\[
H(s) \le H(1) + C (1-s) + (M(1) -M(s)) \le -c_0/2,
\]when $c_0$ is chosen sufficiently large. It follows from (\ref{H'2}) that
\[
H'(s) \ge - c_1 + c_2 H^2(s) - M'(s).
\]This shows
\[
(H(s) + M(s))' \ge - c_3 + c_4 (H(s) + M(s))^2.
\]From here, one immediately deduces
\[
H(1) \ge -A, \qquad A>0.
\]The claim is proven. Thus
\[
\int e^{- \pi |y|^2} \ln G_V(y, s; x, 0) dy \ge -c_0
\]where $|x-y| \le 1$.

Using the reproducing property of $G_V$ and Jensen's inequality,
we have, when $|x-y| \le 1$,
\[
\begin{aligned}
 \ln G_V(x, 2; y, 0)&= \ln \int G_V(x, 2; z, 1) G_V(z, 1; y,
0) dz \\
&\ge \int e^{- \pi |y|^2} \ln G_V(x, 2; z, 1) dz + \int e^{- \pi
|y|^2} \ln G_V(z, 1; y, 0) dz \ge - C.
\end{aligned}
\]This proves the on-diagonal lower bound. The full Gaussian lower
bound now follows from the standard argument in [FS]. \qed

\medskip

{\it Proof of (iii).}

Since,
\[
V = \Delta f - \alpha |\nabla f|^2 - \partial_t f,
\]we have
\[
V + (\alpha-1) |\nabla f|^2 = \Delta f - |\nabla f|^2 -\partial_t f.
\]Let $u = e^{-f}$, by direct calculation,
\[
\Delta u + V u - \partial_t u + (\alpha-1) |\nabla f|^2 u = 0.
\]Hence
\[
u(x, t) = \int G_V(x, t; y, 0) u(x, 0) dy + (\alpha-1) \int^t_0\int
 G_V(x, t;y, s) |\nabla f|^2(y, s) u(y, s) dyds.
\]Since $f \in L^\infty$, we know that $u$ is bounded between two positive constants.
If, by assumption, $G_V$ has a Gaussian lower bound, then, for some
$b>0$, we have
\[
\int^t_0 \int g_b(x, t; y, s) |\nabla f|^2(y, s) dyds \le C
\frac{\sup u}{\inf u}.
\]This completes the proof of part (iii) of Theorem 2.2. \qed

\medskip

{\it Proof of Theorem 2.3.}

(i).  We write
\[
V_1 = 2 ( \Delta f - \partial_t f - 3 | \nabla f|^2), \qquad V_2 = 6
|\nabla f|^2.
\]Let $G_{V_i}$, $i=1, 2$, be the fundamental solution of
$\Delta u + V_i u  - \partial_t u=0$. Since
\[
V = \Delta f - \partial_t f= (V_1/2) + (V_2/2),
\]the Feynman-Kac
formula implies
\[
G_V(x, t; y, s) \le \big{[} G_{V_1}(x, t; y, s) \big{]}^{1/2} \
\big{[} G_{V_2}(x, t; y, s) \big{]}^{1/2}.
\]Observe that
\[
V_1 =  \Delta (2f) - \partial_t(2f)- \frac{3}{2} | \nabla (2f)|^2.
\]Hence, by Theorem 2.2, we know that $G_{V_1}$ has Gaussian upper
bound. Under the smallness assumption on the $L^{p, q}$ norm of
$|\nabla f|^2$ in the theorem, it is well known that $G_{V_2}$ also
has a Gaussian upper bound. Therefore $G_V$ has Gaussian upper
bound.

In order to prove the Gaussian lower bound, we observe that
\[
V=\Delta f - \partial_t f \ge V_3 \equiv \Delta f - \partial_t f - 2
|\nabla f|^2.
\]Under our assumption on the $L^{p, q}$ norm of $|\nabla f|$,
it is straight forward to check that
\[
g_{1/4} \star |\nabla f|^2 \in L^\infty.
\]
Hence Theorem 2.2 (ii) shows $G_{V_3}$ has Gaussian lower bound.
Clearly this Gaussian lower bound of $G_{V_3}$ is also a Gaussian
lower bound of $G_V$ by the maximum principle. This proves part (a).
\medskip

(ii). Clearly we can choose $A_0$ sufficiently small so that all the
following kernels have global Gaussian upper and lower bound:
\begin{equation}
\label{4b}
 G_{2V}, \quad G_{V/2}, \quad G_{2 \Delta f}, \quad
G_{\Delta f/2}.
\end{equation}
The bounds on the first two kernels
follow from part (i). The bounds on the last two kernels follow from
standard theory since $\Delta f = div (\nabla f)$ with $\nabla f$
has a small norm in the suitable $L^{p, q}$ class. (see \cite{LS:1}
e.g.)

Now observe that
\[
-\partial_t f = \Delta f - \partial_t f - \Delta f = V - \Delta f,
\]
\[
\frac{V}{2}= \frac{\Delta f}{2} - \frac{\partial_t f}{2}.
\]By Feynman-Kac formula
\[
G_{(-\partial_t f)} \le \big{(} G_{2V} \big{)}^{1/2} \ \big{(} G_{2
\Delta f} \big{)}^{1/2};
\]
\[
G_{V/2} \le \big{(} G_{2 \Delta f} \big{)}^{1/2} \ \big{(} G_{ -
\partial_t f} \big{)}^{1/2}.
\]Hence (\ref{4b}) show that  $G_{(-\partial_t f)}$ also has global Gaussian upper and
lower bound. Since the setting of the Theorem is invariant under the
reflection $f \to - f$ the result follows. \qed

\medskip
We close this section by giving proofs of the corollaries.

 \medskip

{\it Proof of Corollary 1.}
\medskip

(a). Let $V_k = \min \{ V(x, t), k \}$, $k =1, 2, ...$. Then
(\ref{1.1}) with $V$ replaced by $V_k$ has a unique solution.

Let $J(t) \equiv \int_{\R^n} u^2_k(x, t) dx$.  Then
 \[ J'(t) = 2 \int_{\R^n}
[- \nabla u_k \nabla u_k  + V_k u^2_k ] dx.
\]By our assumption on $V$,
\[
J'(t) \le  2 b J(t)
\]which implies
\[
\int_{\R^n} u^2_k(x, t) dx \le \int_{\R^n} u^2_0(x) dx \ e^{2 b t}.
\]

Therefore if $u_0 \in L^2(D)$, we conclude that $u_k(x, t)$
increases to a finite positive limit $u(x, t)$ as $k \to \infty$,
for all $t$ and for a.e. $x$. Moreover $u(\cdot, t) \in L^2({\bf
R}^n )$.
 We show that the above $u$ is a positive $L^2$ solution to
(\ref{1.1}).
\medskip

Since $u_k$ is a solution to (\ref{1.1}) with $V$ replaced by $V_k$,
for any $\psi \in C^{\infty}_0(\R^n \times (0, T))$, we have
\[
\int ( u_k \psi) |^{t_2}_{t_1} dx -\int^{t_2}_{t_1} \int u_k \psi_t
dxdt - \int^{t_2}_{t_1} \int   u_k \Delta \psi dxdt -
\int^{t_2}_{t_1} \int V_k u_k \psi dxdt =0
\]for
all $t_1, t_2 \in (\delta, t_0)$.

By our assumption $|V_k u_k | \le |V u| \in L^1(\R^n \times (0,
T))$.  Taking $k \to \infty$ and using the dominated convergence
theorem,
 we
obtain
\[
\int ( u \psi) |^{t_2}_{t_1} dx -\int^{t_2}_{t_1} \int u \psi_t dxdt
- \int^{t_2}_{t_1} \int  u \Delta \psi dxdt - \int^{t_2}_{t_1} \int
V u \psi dxdt =0.
\]
This shows that $u$ is a positive solution to (\ref{1.1}). By
Theorem 2.1
\[
V = \Delta f - |\nabla f|^2 - \partial_t f.
\]

(b). Suppose
\[
V = \Delta f - |\nabla f|^2.
\]Due to the $L^2$ convergence, it suffices to
prove that $V_j$ in Definition 2.1 satisfies (\ref{domi}). Let
$\phi$ be a test function, then
\[
\begin{aligned}
\int^\infty_0 \int  V_j \phi^2 dxdt  &= \int^\infty_0
\int  [\Delta f_j - |\nabla f_j|^2] \phi^2 dxdt\\
 &= - 2\int^\infty_0  \int  \nabla  f_j \nabla \phi \ \phi dxdt  -
\int^\infty_0  \int |\nabla f_j|^2 \phi^2 dxdt \\
&\le \int^\infty_0 \int |\nabla \phi|^2 dxdt \phi.
 \end{aligned}
 \]

 \medskip

 (c). The statement is self-evident by part (b) and Theorem 2.1.
\qed

\medskip

{\it Proof of  Corollary 2.}
\medskip

Suppose $V = \Delta f - |\nabla f|^2 + b$. Then, by the same
limiting argument as above, we have
\[
\begin{aligned}
 \int  V \phi^2 dx  &= \lim_{j \to \infty} \int  [\Delta f_j -
  |\nabla f_j|^2] \phi^2 dx + b \int \phi^2 dx \\
 &= \lim_{j \to \infty} \bigg{(}- 2 \int  \nabla  f_j \nabla \phi \ \phi dx  -
 \int |\nabla f_j|^2 \phi^2 dx \bigg{)} + b \int \phi^2 dx .
 \end{aligned}
 \]Therefore
 \[
 \int  V \phi^2 dx \le \int  |\nabla \phi |^2 dx + b \int \phi^2 dx.
 \]Also by part (a) of Corollary 1, (\ref{1.1}) has a positive
  solution when $u_0 \ge 0$.

On the other hand, suppose $V$ satisfies
\[
\int  V \phi^2 dx  \le \int  |\nabla \phi |^2 dx + b \int \phi^2 dx.
\]Write $V_j = \min \{V, j \}$ with $j=1, 2, ...$. Then
\[
\int  (V_j-b) \phi^2 dx  \le \int  |\nabla \phi |^2 dx.
\]Notice that $V_j-b$ is a bounded function. Hence we can apply
Theorem C.8.1 in \cite{Si:1} to conclude that there exists $u_j > 0$
such that
\[
\Delta u_j + (V_j - b) u_j =0.
\]Writing $f_j = - \ln u_j$, we have
\[
V_j = \Delta f_j - |\nabla f_j|^2 + b.
\]By definition, this means
\[
V =  \Delta f - |\nabla f|^2 + b.
\] \qed

\section{heat bounded functions and the heat equation}

Here we introduce another class of singular functions that has its
origin in the Kato type class. As mentioned in the introduction, a
function is in a Kato type class if the convolution of the absolute
value of the function and the fundamental solution of Laplace or the
heat equation is bounded. Here we generalize this notion by a simple
but key stroke, i.e., we delete the absolute value sign on the
function in the definition of the Kato class. More precisely, we
have

\medskip

\noindent {\bf Definition 3.1}. Let $f=f(x, t)$ be a local $L^1$
function in space time and $G_0$ be the standard Gaussian in ${\bf
R}^n$. We say that $f$ is {\it heat bounded} in a domain $\Omega
\subset {\bf R}^n \times {\bf R}^1$ if
\[
G_0 \star f(x, t) \equiv \int^t_0\int_{{\bf R}^n} G_0(x, t; y, s)
f(y, s) dyds
\]is a bounded function in $\Omega$.

We say that $f$ is {\it almost heat bounded} in a domain $\Omega
\subset {\bf R}^n \times {\bf R}^1$ if
\[
 G_0 \star f(x, t) \in L^p(\Omega)
\]for all $p>1$.

\medskip

{\bf Example}. The function $V(x) = a \frac{\chi_{B(0, 1)}}{|x|^2}$
is not heat bounded but is almost heat bounded in ${\bf R}^n$. Here
$a$ is a nonzero constant.

\medskip

In the next two propositions, we provide a comparison between the
heat bounded class and more familiar classes of functions.

\medskip

\begin{proposition}
 \label{p3.1}
 Suppose, in the distribution sense,
$V =
\partial^2_{ij} f$ with $f \in \cap_{p>1} L^p({\bf R}^n \times (0, T))$. Then
$V$ is almost heat bounded in ${\bf R}^n \times (0, T)$. \proof
\end{proposition}

Let $G_0$ be the free heat kernel in ${\bf R}^n \times (0, \infty)$.
By the assumption on $f$, the function $u=u(x, t)$, defined by
\[
u(x, t) = \int^t_0\int_{{\bf R}^n} G_0(x, t; y) f(y, s) dyds
\] is a solution to the equation
\[
\begin{cases}  \Delta u(x, t) - u_t(x, t) =  - f(x, t), \quad x \in {\bf R}^n,
t>0;\\
u(x, 0) = 0.
\end{cases}
\]By the parabolic version of the Calderon-Zygmond inequality (see
[Lieb] e.g., we know that
\[
u \in W^{2, p}({\bf R}^n), \quad \forall p>1.
\]Hence
\[
\partial^2_{ij} \int^t_0\int_{{\bf R}^n} G_0(x, t; y) f(y, s) dyds
\in L^p({\bf R}^n), \quad \forall p>1. \qed
\]
\medskip

\begin{proposition}
\label{3.2} Suppose,   $0 \le V \in L^1_{loc}$ is form bounded in $D
\times [0, T]$. i.e.
\[
\int^T_0 \int_D  V \phi^2 \le b_1 \int^T_0 \int_D |\nabla \phi|^2 +
b_ 2 \int^T_0 \int_D \phi^2.
\]for all smooth, compactly supported function $\phi \in D \times [0, T]
\subset \R^n \times [0, T]$. Then $V$ is almost heat bounded in $D
\times [0, T]$.
 \proof
\end{proposition}

We will only consider the case when $D= {\bf R}^n$. The other cases
follow from the full space case by a standard comparison method.

Since one can consider $c V$ with $c$ sufficiently small otherwise,
we can choose the constant $b_1$ in the definition of form
boundedness to be $1/2$, i.e. we assume that
\[
\int^T_0 \int_D  V \phi^2 \le \frac{1}{2}\int^T_0 \int_D |\nabla
\phi|^2 + b_ 2 \int^T_0 \int_D \phi^2.
\]for all smooth, compactly supported function $\phi$.

Let $u_k$ be the solution of
\begin{equation}
\label{3.1}
\begin{cases}
 \Delta u_k + V_k u_k - (u_k)_t=0, \quad \text{in} \quad {\bf R}^n \times
(0, \infty), \quad V \in L^2_{loc}({\bf R}^n \times
(0, \infty))\\
u_k(x, 0) = u_0(x)>0, \quad x \in {\bf R}^n, \quad u_0 \in L^2({\bf
R}^n ).
\end{cases}
\end{equation}Here $V_k$ is the truncated potential $V_k = \min \{ V, k \}$
with $k$ being positive integers. Clearly $V_k \le V_{k+1}$.
\medskip

 We show that $u_k$ converge pointwise to a locally integrable
function.

Let $J(t) \equiv \int_D u^2_k(x, t) dx$.  Then
 \[ J'(t) = 2 \int_D
[- \nabla u_k \nabla u_k  + V_k u^2_k ] dx.
\]By our assumption on $V$,
\[
J'(t) \le  2 b J(t)
\]which implies
\[
\int_D u^2_k(x, t) dx \le \int_D u^2_0(x) dx \ e^{2 b t}.
\]

Therefore if $u_0 \in L^2(D)$, we conclude that $u_k(x, t)$
increases to a finite positive limit $u(x, t)$ as $k \to \infty$,
for all $t$ and for a.e. $x$. Moreover $u(\cdot, t) \in L^2({\bf
R}^n )$.

Write $w_k = \log u_k$. From (\ref{3.1}), one deduces
\[
\Delta w_k + |\nabla w_k|^2 + V_k -(w_k)_t =0.
\]Therefore
\[
\begin{aligned}
w_k(x, t)& = \int_D G_0(x, t; y, 0) w_k(x, 0) dy \\
&\qquad + \int^t_0\int_D G_0(x, t; y, s) |\nabla w_k(x, s)|^2 dyds +
\int^t_0\int_D G_0(x, t; y, s) V_k(y, s) dyds.
\end{aligned}
\]Therefore
\[
\int^t_0\int_D G_0(x, t; y, s) V_k(y, s) dyds \le w_k(x, t) - \int_D
G_0(x, t; y, 0) w_k(x, 0) dy.
\]By the monotone convergence theorem
\[
\begin{aligned}
\int^t_0\int_D G_0(x, t; y, s) V(y, s) dyds &\le w(x, t) - \int_D
G_0(x, t; y, 0) w(y, 0) dy \\
&\le \log (1 + u(x, t)) -  \int_D G_0(x, t; y, 0) \log u_0(y) dy.
\end{aligned}
\]Now we take
\[
u_0(x) = \frac{1}{1+|x|^n}.
\]Then, since  $u_0 \in L^2({\bf R}^n)$, we have
\[
u(\cdot, t) \in L^2({\bf R}^n).
\]By Jensen's inequality,
\[
\log (1 + u(\cdot, t)) \in L^p({\bf R}^n), \qquad \forall p>1.
\]It is also clear that
\[
 \int_D G_0(\cdot , t; y,
0) \log (1+|y|^n) dy \in L^p_{loc} ({\bf R}^n), \qquad \forall p>1.
\]The result follows. \qed

\medskip

\medskip

\section{applications to the Navier-Stokes equation}
In this section, we establish a new a priori estimate for a certain
quantity involving the velocity and vorticity of the $3$ dimensional
Navier-Stokes equation.

\[
\aligned
 &u_t - \Delta u(x, t) + u \cdot \nabla u(x, t) + \nabla p = 0,\\
 &\nabla \cdot u = 0,\\
 &u(x, 0) = u_0(x)
\endaligned
 \leqno (4.1)
\]
for $(x, t) \in {\mathbb R}^3 \times (0, \infty)$, where $\Delta$ is
the standard Laplacian, a vector field $u$ represents the velocity
of the fluid, and a scalar field $p$ the pressure. (The viscosity is
normalized, $\nu=1$.)

There has been an extensive and rapidly growing literature on the
equation which is impossible to quote extensively here. Let us just
mention that weak solutions are known to exist due to the seminal
work of Leray \cite{L:1}. However it is not known if the weak
solution is smooth everywhere. Several sufficient conditions
implying smoothness of weak solutions have been made. See for
example \cite{P:1} and \cite{S:1}. In these two papers, it was shown
that if the velocity $u$ is in $L^{p, q}$ class with $\frac{3}{p} +
\frac{2}{q} < 1$, then $u$ is actually smooth. For more sufficiency
results in various other spaces we refer the reader to the more
recent survey paper \cite{Ca:1}. However it is only known that $u
\in L^{10/3, 10/3}$. Therefore there is a gap in between the a
priori estimate and the sufficiency condition.

What we will prove here is a different sufficiency condition and a
priori estimate using the heat bounded and almost heat bounded
potentials defined in the previous section. There is still a gap
between the two conditions. However the gap seems logarithmic. More
precisely, we have

\medskip

\begin{theorem}
\label{iff} Let $u$ be a Leray-Hopf solution of the Navier-Stokes
equation, which is classical in ${\bf R}^3 \times (0, T)$.  Let $w$
be the vorticity $\nabla \times u$. Define the quantity
\[
{\bf Q} = {\bf Q}(x, t) \equiv \frac{curl( u \times w) \cdot w + 2 |
\nabla \sqrt{|w|^2+1} |^2 - | \nabla w|^2}{| w |^2+1}(x, t).
\]

Then the following statements hold for any $\delta \in (0, T)$.

(1). The quantity ${\bf Q}$ is almost heat bounded in $ \left( {\bf
R}^3 \times (\delta, T]\right) \cap \{ |w| \ge 1 \}.$

(2). $u$ is a classical solution of the Navier-Stokes equation in
${\bf R}^3 \times (\delta, T]$ if and only if ${\bf Q}$ is heat
bounded in $ \left( {\bf R}^3 \times (\delta, T]\right) \cap \{ |w|
\ge 1 \}.$
\end{theorem}
\medskip

{\it Remark 4.1.}

The quantity $Q$ is well defined since we assume that $u$ is smooth
for $t \in (0, T)$. The first term in $Q$ is essentially the vortex
stretching factor which is the hardest to control. The point of the
theorem is that if there is blow up at time $T$, then the blow up
just happens barely.

\medskip

{\it Proof of Theorem 4.1.}

We will just prove (1) since (2) is self-evident afterward.

We divide the proof into three steps.
\medskip

{\bf Step 1.} rewriting the vortex equation in the log form.
\medskip

Let $w=w(x, t)$ be the vortex. It is well known that $|w|^2$
satisfies the following scalar heat equation with lower order terms
\[
\Delta |w|^2 - u \cdot \nabla |w|^2 + 2 \alpha |w|^2 - 2 |\nabla
w|^2 -(|w|^2)_t = 0. \leqno(4.2)
\]Here $\alpha$ is the vortex stretching potential given by (c.f.
[Co])
\[
\aligned
 \alpha(x, t) &= \frac{3}{4 \pi} P.V. \int_{{\mathbb R}^3}
D[\widetilde{y}, \widetilde{\w}(x+y), \widetilde{\w}(x)] \ |\w(x+y,
t)| \frac{dy}{|y|^3}\\
& =  \frac{ w \nabla u \cdot w}{| w |^2}.
 \endaligned\leqno(4.3)
\]

 A straightforward computation from (4.2) shows
\[
\Delta \ln (|w|^2 +1) - u \cdot \nabla  \ln (|w|^2 +1) + 2
\frac{\alpha |w|^2}{|w|^2+1} - 2 \frac{|\nabla w|^2}{|w|^2+1} +
\frac{|(\nabla |w|^2)|^2}{(|w|^2+1)^2} - \partial_t (\ln (|w|^2+1))
= 0.
\]Write $f = \frac{1}{2} \ln (|w|^2+1)$. We deduce
\[
\Delta f - u \cdot \nabla  f + \frac{\alpha |w|^2}{|w|^2+1} + 2
|\nabla f|^2 - \frac{|\nabla w|^2}{|w|^2+1}
 -f_t =
0. \leqno(4.4)
\]
\medskip

{\bf Step 2.} a representation formula.

\medskip

By our assumption, for $t \in (0, T)$, $u$ and $w$ are classical
functions and $f$ vanishes near infinity.  This shows,
\[
\aligned
 f(x, t) &= \int G_0(x, t; y, s) f_0(y) dy \\
 & \qquad + \int^t_0 \int
G_0(x, t; y, s) \left[ \frac{\alpha |w|^2}{|w|^2+1} - u \cdot \nabla
f + 2 |\nabla f|^2 - \frac{|\nabla w|^2}{|w|^2+1}
 \right] (y, s) dy ds.
\endaligned
\leqno(4.5)
\]
\medskip

{\bf Step 3.}  Apply Jensen's inequality.

\medskip

For convenience, we write
\[
Q \equiv \frac{\alpha |w|^2}{|w|^2+1} - u \cdot \nabla f + 2 |\nabla
f|^2 - \frac{|\nabla w|^2}{|w|^2+1} \leqno(4.6)
\]It is clear that
\[
\aligned
 Q & = \frac{ w \nabla u \cdot w}{| w |^2+1} - \frac{1}{2} u_j
\partial_j \ln (|w|^2+1) + 2 |\nabla f|^2 - \frac{|\nabla
w|^2}{|w|^2+1}\\
& = \frac{ u \nabla w \cdot w}{| w |^2} - \frac{u_i
\partial_j w_i w_i}{|w|^2+1} + 2 |\nabla f|^2 -  \frac{|\nabla
w|^2}{|w|^2+1} \\
&= \frac{ w \nabla u \cdot w - u \nabla w \cdot w + 2 | \nabla
\sqrt{|w|^2+1} |^2 - |\nabla w|^2 }{| w |^2+1}.
\endaligned
\]Following the well known vector identity, we have
\[
Q = \frac{curl( u \times w) \cdot w + 2 | \nabla \sqrt{|w|^2+1} |^2
- |\nabla w|^2}{| w |^2+1}. \leqno(4.7)
\]

It is well known that $w \in L^2({\bf R}^3 \times \R^+)$. Using
Jensen's inequality, it is easy to show that $f(\cdot, t) \in L^p$
for any $p>1$, in the region where $|w| \ge 1$. Hence the quantity
$Q$ is almost heat bounded.
\qed

\bigskip

\noindent e-mail: qizhang@math.ucr.edu

\enddocument